\documentclass[11pt,letterpaper]{amsart}
\usepackage{hyphenat}
\usepackage[T1]{fontenc}
\usepackage{lmodern}
\usepackage{amssymb,mathrsfs}
\usepackage[margin=1.1in]{geometry}
\usepackage{microtype}
\usepackage{needspace}
\usepackage{constants}
\usepackage{bm} 
\usepackage[colorlinks=true,linkcolor=blue,citecolor=black,urlcolor=blue,
linktocpage=true]{hyperref}
\hypersetup{pdftitle={Supercritical level-set percolation on finite graphs},
	pdfsubject={Probability theory},
	pdfkeywords={Gaussian free field, percolation, spectral gap, giant cluster}}

\newconstantfamily{c}{symbol=c}    

\newtheorem{theorem}{Theorem}[section]
\newtheorem{proposition}[theorem]{Proposition}
\newtheorem{lemma}[theorem]{Lemma}
\newtheorem{corollary}[theorem]{Corollary}
\numberwithin{equation}{section}
\renewcommand{\P}{\mathbb P}

\begin{document}
	\title[Supercritical level-set percolation]{Supercritical level-set percolation on finite graphs}
    \author[S.~Goswami]{Subhajit Goswami}
      \address{School of Mathematics, Tata Institute of Fundamental Research, 1, Homi Bhabha Road, Colaba, Mumbai 400005, India.}
      \curraddr{}
       \email{goswami@math.tifr.res.in}
       \thanks{}

     \author[D.~Pal]{Dipranjan Pal} 
     \address{Statistics and Mathematics Unit (SMU), Indian Statistical Institute, Kolkata, 203 B.T. Road, Kolkata 700108, India.}
     \curraddr{}
     \email{dipranjan\_r@isical.ac.in}
     \thanks{}
    
	\date{}
	\begin{abstract}
		We continue the study of bond percolation on finite graphs induced by the level-sets of the zero-average Gaussian free field on the associated metric graph. In our previous work 
        \cite{goswami2026crifinite}, we analyzed the critical and near-critical behavior of this model on uniform expanders, focusing on the volume of the largest cluster. In the present paper, we investigate the size of the second largest cluster in the supercritical regime --- the phase where one expects a unique giant component. We prove that the second largest cluster is at most of 
        logarithmic order with high probability. The estimates hold uniformly over deterministic graphs with a bounded maximum degree and a spectral gap bounded away from zero. In particular, we do not require any assumption regarding regularity, girth or local geometry. This picture is substantially different from Bernoulli percolation, where the second largest supercritical cluster can have polynomial volume under identical hypotheses. The proof relies on expansion properties of the underlying graph and the Dirichlet martingale framework developed in \cite{goswami2026crifinite}. 
	\end{abstract}
	\maketitle
	\section{Introduction}\label{sec:intro}
\subsection{Background and motivation}\label{sec:back}
Percolation on finite graphs features many interesting 
phenomena that are related to, yet fundamentally distinct from 
percolation on infinite graphs. A central question concerns the 
emergence and uniqueness of a giant component --- a cluster 
comprising a positive fraction of the vertices. While the answer to this question is known classically for the Erd\H{o}s-R\'enyi model \cite{MR148055}, cluster behavior across general families of graphs can depend substantially on the underlying geometry. In 
the present article, we revisit this question for the bond percolation model introduced in
\cite{goswami2026crifinite}, obtained from the level-sets of the zero-average Gaussian free field on finite metric graphs. 
Our interest lies in results that hold uniformly over bounded-degree expander families, without any further assumptions on 
their local geometry.

For Bernoulli bond percolation, Alon, Benjamini and Stacey
\cite{MR2073175} proved that a family of expanders with 
uniformly bounded degrees has, with high probability, at most 
one cluster of any given linear size, uniformly in the percolation parameter. A different line of work concerns
finite vertex-transitive graphs. Easo and Hutchcroft \cite{MR4801596}
established uniqueness in the supercritical regime for every sparse family
of such graphs, thereby resolving a conjecture of Benjamini \cite{benjamini2001percolation}, and
characterized the obstruction to uniqueness in the dense case.

A logarithmic bound on the size of the second largest cluster 
in the supercritical regime --- known to hold for Erd\H{o}s-R\'enyi and random regular graphs \cite{MR756039, MR2435852, MR2583058} --- is a substantially stronger property than the 
mere uniqueness of the giant component. For Bernoulli 
percolation on bounded-degree expander graphs $G = (V, E)$, a polynomial upper bound $|V|^\omega$ with $\omega < 1$ on the second largest component was established in \cite{MR2073175}. 
Later, Krivelevich, Lubetzky and Sudakov \cite{MR4101348} 
showed that polynomial sizes can in fact occur: even among regular 
expanders of arbitrarily large girth, the second largest 
supercritical cluster can have size at least $|V|^\beta$ for 
any prescribed $\beta\in(0,1)$, provided the expansion constant is chosen suitably. More recently, Diskin and Krivelevich \cite{MR4922676} proved a logarithmic bound for
regular expander graphs of fixed degree under an additional 
sparsity condition on sets of polylogarithmic size. Their companion 
work \cite{diskin2025componentslargesmalli} treats growing 
degrees under suitable global and small-set expansion 
assumptions. These results shed light on the additional 
geometric information that can yield logarithmic bounds.

Beyond Bernoulli percolation, level-sets of variants of the 
Gaussian free field provide a natural class of percolation 
models with long-range correlations. This model has been studied extensively in 
recent years on $\mathbb Z^d$ ($d \ge 3$) and other infinite 
transient graphs; see \cite{RodriguezSznitman13, MR4112719, gosrodsev2021radius, MR4557402, DCGRS20, drewitz2026fractal} 
and references therein. For finite graphs, where the 
unconstrained Green function is infinite, zero-mean normalization yields a canonical construction: the so-called zero-average Gaussian free field (see \S\ref{subsec:prelim_field} for a formal definition). Ab\"acherli \cite{abacherli2018local} studied its level-set percolation on large discrete tori. Ab\"acherli and \v{C}ern\'y \cite{MR4115734,MR4169171} subsequently considered regular trees
and a class of locally tree-like regular expanders. In the 
latter setting, they identified the threshold separating 
logarithmic subcritical clusters from a supercritical phase in which a positive fraction of the vertices belong to mesoscopic 
clusters. \v{C}ern\'y \cite{MR4642820} proved that these
clusters give rise to a unique giant throughout the supercritical phase without any quantitative bound for the 
volume of the second largest cluster. The class of graphs considered there includes large-girth regular expanders and typical
realizations of random regular graphs. For random 
regular graphs, Conchon-Kerjan \cite{MR4552959} obtained a more detailed description under the {\em annealed} law for the joint 
realization of the graph and the field, including uniqueness of 
the giant and a second largest cluster of order $\log |V|$. In 
the same paper, the author conjectured that the logarithmic 
order of the non-giant cluster for the level-set percolation on deterministic regular graphs would require 
more than mere expansion properties, similar to the case of 
Bernoulli percolation (see Section~1.4). Our results in the 
current article show that, at least for the level-sets of the 
associated metric graph, this is {\em not} the case.

The metric graph formulation of free-field level-set percolation originates in the work of Lupu \cite{MR3502602}. In \cite{goswami2026crifinite}, the authors of the current 
paper introduced the corresponding zero-average model on 
arbitrary finite graphs and established various tail bounds on 
the tail probabilities for the size of the largest cluster 
depending only on the number of vertices, the
maximum degree and the spectral gap. For uniformly bounded degrees and
uniformly positive spectral gaps, these bounds identify the critical level
as $0$, give clusters of order $|V|^{2/3}$ in a window of width
$|V|^{-1/3}$, and yield a giant at every fixed negative level. In the parlance of random graphs, these properties 
characterize the {\em mean-field} (near-)critical behavior. The zero-average capacity and the associated 
Dirichlet martingale framework developed in \cite{goswami2026crifinite} allow one to 
work with global properties of the graph, without requiring 
much control on the local geometry.

The purpose of the present article is to describe the supercritical phase
more precisely by proving uniqueness of the giant together with a
logarithmic upper bound for the second largest cluster. The distinction from
the preceding results is the generality of the underlying deterministic
graphs: the bounds require only a uniform degree bound and a uniform lower
bound on the spectral gap. They impose no regularity, transitivity,
large-girth or local convergence assumption. In particular, the
logarithmic bound holds under hypotheses that permit polynomially large
second largest clusters for Bernoulli percolation. This feature makes the metric
graph zero-average model a setting in which one can study the supercritical
phase uniformly over a broad class of graphs through their global
isoperimetric properties.

\subsection{Main result}\label{subsec:results}
Let $\mathcal G=(V,E)$ be a finite, simple, connected graph. 
Its zero-average Gaussian free field $\varphi$ is the centered
Gaussian field with covariance given by the zero-average Green 
function; the normalization and its extension to the metric 
graph are recalled in Section~\ref{sec:preliminaries}. Here, as throughout the paper, we leave the dependence on 
$\mathcal G$ implicit in all our notations. Given the field on $V$, an edge $\{x,y\}\in E$ is declared {\em open} at level $h\in\mathbb R$, independently of the
other edges, with conditional probability
\begin{equation}\label{eq:opening-rule}
 1-\exp\{-2(\varphi(x)-h)_+(\varphi(y)-h)_+\},
 \qquad a_+=\max\{a,0\}.
\end{equation}
This is precisely the probability that the Brownian bridge 
joining its endpoint values on the corresponding unit interval 
stays above $h$. We write $E^{\ge h}$ for the open edges and
$\mathcal G^{\ge h}=(V,E^{\ge h})$. We denote its cluster containing $x$ by $\mathscr C_x^{\ge h}$, and its largest and second largest clusters
by $\mathscr C_{\max}^{\ge h}$ and $\mathscr C_{\sec}^{\ge h}$,
respectively. Ties are resolved using some fixed, deterministic ordering. If there is only one cluster, we put 
$|\mathscr C_{\sec}^{\ge h}| = 0$. 

Let $D_{\max}$ and $\lambda_s$ denote the maximum
degree and the normalized spectral gap of $\mathcal G$, 
respectively; see \eqref{def:spectral_gap}. All relevant random objects --- the field $\varphi$ and the variables $\{1_{\{e \in E^{\ge h}\}} : e \in E, h \in \mathbb R\}$ --- are 
defined on a common probability space $(\Omega, \mathcal F, \P)$. 

\Needspace{8\baselineskip}
\begin{theorem}\label{thm:uniq_giant}
Let $D \ge 3$, $\lambda > 0$ and $h < 0$. There exist
$c = c(h, D, \lambda) > 0$ and $C = C(h, D, \lambda) < \infty$ 
such that every finite, simple, connected graph $\mathcal G=(V,E)$ with $D_{\max} \le D$ and $\lambda_s \ge \lambda$ satisfies
\begin{equation}\label{eq:uniq_giant}
 \P \big[
 |\mathscr C_{\max}^{\ge h}| < c |V|
 \ \text{or}\ 
 |\mathscr C_{\sec}^{\ge h}|>C\log |V|\big]
 \le C|V|^{-1/15}.
\end{equation}
\end{theorem}

Thus, for every family satisfying these two uniform graph bounds and
with $|V|\to\infty$, there is a unique giant cluster at every fixed
$h<0$, and all other clusters have volume at most $C\log|V|$, with
probability tending to one. The exponent in \eqref{eq:uniq_giant} comes from %the giant-existence estimate of
the probability for the existence of a giant cluster as stated 
in \cite[Theorem 1.2]{goswami2026crifinite}. 

%The two new estimates below, uniqueness at a fixed linear scale and absence of clusters between
%logarithmic and linear scales, have error at most $C|V|^{-2}$.

\subsection{Discussion of the proof}\label{subsec:discussion}
There are two distinct issues. Large clusters must be joined 
by a small sprinkling, and one must exclude every cluster whose
volume lies between logarithmic and linear order. Both steps 
use the global expansion of the underlying graph, but in 
different ways.

\smallskip

For the first step, we extract a small independent Gaussian 
component from the field. Conditionally on the remaining Gaussian field, this perturbation has independent Gaussian variables at the vertices with a common {\em centering} term. 
A similar decomposition was also used in the work \cite{MR4642820} by \v{C}ern\'y. Expansion property of $\mathcal G$ --- as captured by the lower bound on $\lambda_s$ --- supplies linearly many vertex-disjoint paths of bounded 
length between any two linear size clusters. One can then 
leverage the independence of the extracted field and 
properties of Brownian bridges, together with a delicate 
control on the effect of the centering term on the clusters, to deduce exponential in length lower bounds on the 
conditional probabilities of opening any such path at a 
slightly lower level. The linear number of paths 
then allow us to connect every pair of large clusters at this 
new level with high probability.

For the second step, we use the zero-average capacity and the Dirichlet martingale introduced in \cite{goswami2026crifinite}. The spectral gap gives a linear-in-volume lower bound for the zero-
average capacity proportional. At the time when the first cluster exploration is complete, its boundary
values are at most the level. A Gaussian estimate makes the bulk part of the
martingale small for all sets below a suitably chosen linear volume. For any explored cluster of such a volume, the martingale must therefore lie below a negative linear function of its capacity. The capacity is also the quadratic variation of the martingale which, by Dubin-Schwarz, corresponds to its 
`Brownian time'. The event that a Brownian motion crosses this line after time proportional to the cluster volume has an exponentially small probability. A union bound over the 
starting vertex yields the logarithmic cutoff. This
uses only the exploration time and requires no local 
approximation of the graph.

\subsection{Comparison with earlier methods}\label{subsec:methods-comparison}
	We now explain how the preceding arguments relate to earlier approaches.
	For the discrete vertex level-sets on regular expanders,
	\cite{MR4169171} and \cite[Sections 4--7]{MR4642820} compare the field
	locally with the free field on the regular tree. This comparison supplies
	mesoscopic clusters to which sprinkling can be applied. In the random
	regular setting, Conchon-Kerjan \cite[Sections 1.5 and 3--7]{MR4552959}
	works under the joint law of the graph and the field (i.e., the annealed law), using random-graph
	pairings and a recursive approximation of the field. He couples the field on the explored sets with its tree 
    counterpart. In our setting the graph
	is fixed and its neighborhoods need not resemble trees. We 
    use the existence of giant clusters as guaranteed by the results in \cite{goswami2026crifinite} and control
	smaller clusters through zero-average capacity, which is 
    bounded below in terms of volume by the spectral gap.
	
	The use of expansion via Menger's theorem to join large clusters by sprinkling has a classical precedent in Bernoulli percolation; see, for example, \cite[Proposition 3.1]{MR2073175}. The extraction of an independent Gaussian 
    part is similar in spirit to the approach of \cite[Sections 3 and 9]{MR4642820} although our proof 
    strategies are different.
	%In that work the decomposition is available on general finite graphs; tree comparison enters the construction of 
    %mesoscopic clusters.
    
	 We use the Dirichlet martingale framework developed in \cite{gosrodsev2021radius}, inspired by the work \cite{MR4112719}, to exclude intermediate
	clusters without examining their local geometry. For comparison to our two-step approach outlined in the previous subsection which is inspired by classical arguments for Erd\H{o}s-R\'enyi random graphs, the works of Diskin and Krivelevich 
    \cite{diskin2025componentslargesmalli,MR4922676} in the Bernoulli case also separate
	the joining of large clusters from the exclusion of intermediate ones.
	In the growing-degree case, the latter uses spanning-tree counts and
	closed-boundary probabilities
	\cite[Lemma 3.4]{diskin2025componentslargesmalli}; at fixed degree it
	uses a modified breadth-first exploration and a binomial deviation
	estimate \cite[Lemma 3.1]{MR4922676}. Their
	stronger assumptions on small sets provide the quantitative control
	needed in these estimates. In our bounded-degree argument this role is played by the capacity estimate. Expansion alone already suffices for Bernoulli uniqueness on the linear scale
	\cite[Section 2]{MR2073175}; the distinction is the logarithmic bound for the remaining clusters.

\medskip

	In Section~\ref{sec:preliminaries} we recall the field, zero-average
	capacity and martingale functional. Section~\ref{sec:auxiliary} establishes some auxiliary estimates required for the proof of our main result. Section~\ref{sec:supcri} is devoted to the proof of Theorem~\ref{thm:uniq_giant}.

\medskip

\noindent\emph{Notations and conventions.} We denote by $c, c', C, C', \ldots$ etc. finite, positive constants 
which are allowed to vary from place to place. All the constants implicitly depend on $D$ and $\lambda$. Their 
dependence on any parameters will be made explicit. Numbered constants $c_1, C_1,\dots$ stay fixed after their first appearance 
within the text. We fix the graph $\mathcal G = (V, E)$ and usually 
omit it from the notation. We use $|\bm{\cdot}|$ to denote the cardinality of  $\bm{\cdot}$ when $\bm{\cdot}$ is a set and 
the absolute value of $\bm{\cdot}$ when it is a real number.

\section{Preliminaries}\label{sec:preliminaries}
We recall the relevant constructions and the potential-theoretic 
objects from \cite[Sections 2--4]{goswami2026crifinite}.
\subsection{The graph, cable system and the zero-average field}\label{subsec:prelim_field}
The {\em Dirichlet form} on $\mathcal G$ is given by
\begin{equation}\label{eq:discrete-energy}
 \mathcal E(f,g) \stackrel{\rm def.}{=} \sum_{\{x,y\}\in E}(f(x)-f(y))(g(x)-g(y)),
\end{equation}
where $f,g \in \mathbb R^V$. The (normalized) {\em spectral gap} $\lambda_s$ of $\mathcal G$ can be defined as
\begin{equation}\label{def:spectral_gap}
\lambda_s = \inf\Big\{\mathcal E(f,f): \sum_{x\in V} {\rm deg}(x) f^2(x) =1, \, \sum_{x \in V} {\rm deg}(x) f(x) = 0 \Big\}
\end{equation}
where ${\rm deg}(x)$ denotes the degree of the vertex $x$ 
in the graph $\mathcal G$. As required by 
Theorem~\ref{thm:uniq_giant}, we assume throughout that $\max_{x \in V} {\rm deg}(x) \le D$ and $\lambda_s \ge \lambda>0$. In particular, by Cheeger’s inequality (see,~e.g.,~\cite[Theorem~6.15]{MR3616205}), we have
\begin{equation}\label{eq:edge-expansion}
 |\partial_e S|\ge c \min\{|S|, |V \setminus S|\},\qquad S \subset V,
\end{equation}
where $\partial_e S$ is the {\em edge boundary} of $S$, i.e., the 
set of edges with exactly one endpoint in $S$.

\smallskip

% For completeness, the covariance of the field can be specified 
% without a choice of random-walk time scale. Let $\Delta$ denote 
% the graph {\em Laplacian} associated with \eqref{eq:discrete-energy}, and let $L^\dagger$ be its inverse on the
% Euclidean orthogonal complement of the constant functions, extended as zero
% on the constants. Set
% \begin{equation}\label{eq:projection}
%  \pi_x=\frac{{\rm deg}(x)}{q},\qquad
%  \Pi=I-\mathbf1\pi^{\mathsf T},\qquad
%  G=\Pi L^\dagger\Pi^{\mathsf T}.
% \end{equation}
The {\em zero-average Gaussian free field} (GFF) $\varphi=(\varphi(x))_{x\in V}$ is the centered Gaussian process with 
covariances given by the {\em zero-average Green function} $g(\cdot, \cdot)$ on $\mathcal G$ (see \cite[\S2.4]{goswami2026crifinite} and references therein). Equivalently, it is the Gaussian
field associated with the form \eqref{eq:discrete-energy} on the Hilbert space $\{f \in \mathbb R^V: \sum_{x \in V} {\rm deg}(x) f(x)=0\}$. We recall from \cite[Lemma~2.1]{goswami2026crifinite} that
\begin{equation}\label{eq:weighted-sum-variance}
 {\rm Var}\Big[\sum_{x \in S} {\rm deg}(x) \varphi(x)\Big]\le C|S|, \qquad S\subset V.
\end{equation}
% In particular,
% \begin{equation}\label{eq:weighted-zero-average}
%  \sum_{x\in V} {\rm deg}(x) \varphi(x) = 0, \quad \text{$\P$-almost surely}.
% \end{equation}
% The spectral-gap inequality also gives, for every $v\in\mathbb R^V$,
% \begin{equation}\label{eq:covariance-upper}
%  v^{\mathsf T}Gv\le\lambda^{-1}\sum_{x\in V}v_x^2.
% \end{equation}
% Indeed, the variance on the left is the supremum of
% $(\sum_xv_xf(x))^2/\mathcal E(f,f)$ over nonzero functions of
% degree-weighted average zero, and \eqref{eq:spectral-gap} bounds the
% denominator below by $\lambda\sum_xf(x)^2$.

\smallskip

The {\em metric graph} or the {\em cable system} 
$\widetilde{\mathcal G}$ is obtained by replacing each edge $e=\{x,y\}$ by an open interval $I_e$ of unit length, with endpoints 
identified with $x$ and $y$. We equip it with the path metric, in 
which the distance between two points is the minimum length of a 
path joining them. We write $B(x,r)$ for the open ball of radius $r$ 
centered at $x$ in this metric, and $\partial K$ for the topological 
boundary of a set $K$. Conditionally on the vertex values $(\varphi(x) : x \in V)$, we extend $\varphi$ independently along the different edges by standard Brownian bridges joining their endpoint values. We use the same notation $\varphi$ for the resulting continuous field and $\P$ for the underlying probability measure. For any $h \in \mathbb R$, let 
\begin{equation}\label{eq:cable-level-set}
 \widetilde{\mathcal G}^{\ge h} \stackrel{{\rm def.}}{=} \{x \in\widetilde{\mathcal G}:\varphi(x)\ge h\} \: \mbox{ and }
 \: E^{\ge h} \stackrel{{\rm def.}}{=} \{e\in E:\overline I_e\subset \widetilde{\mathcal G}^{\ge h}\}
\end{equation}
denote the {\em level-set above} height $h$ and the corresponding 
set of {\em open} edges respectively (see \S\ref{subsec:results}). 
Here $\overline I_e$ denotes the topological closure of $I_e$. The 
reflection principle then gives the conditional opening probability 
in \eqref{eq:opening-rule}. See \cite[Section~2]{goswami2026crifinite} for more details. %The same bridges couple these graphs monotonically for all levels $h$.

For any $x \in \widetilde{\mathcal G}$ and $S \subset \widetilde{\mathcal G}$, let $\mathcal C_x(S)$ denote the {\em 
(connected) component} of $x$ inside $S$ (as subspaces of 
$\widetilde{\mathcal G}$). Further for any $S' \subset S$, let $\mathcal C_{S'}(S) \stackrel{{\rm def.}}{=} \cup_{x \in S'} \, \mathcal C_x(S)$. %For $S,A\subset\widetilde{\mathcal G}$, let $C_A(S)$ denote the union of the connected components of $S$ that intersect $A$. We write $C_A^{\ge u}(S)=C_A(S\cap\widetilde{\mathcal G}^{\ge u})$ and $|K|_V=|K\cap V|$. 
As in \cite{goswami2026crifinite}, we use \emph{cluster} for a connected component of the random graph $\mathcal G^{\ge h} = (V,E^{\ge h})$ and reserve 
the term \emph{component} for maximal connected subsets 
of $\widetilde{\mathcal G}$. %In particular, a vertex below $h$ is an isolated cluster in $(V,E^{\ge h})$, while its component in the metric level-set is empty. This distinction has no effect on statements about clusters with at least two vertices.

\subsection{Zero-average capacity and the Dirichlet martingale}\label{subsec:zvgcap}
Let $\mathcal K_{<\infty}$ denote the collection of nonempty compact 
subsets $K$ of $\widetilde{\mathcal G}$ having finitely many 
components and satisfying $V\setminus K\ne\varnothing$. On the metric graph, the Dirichlet form is
\begin{equation}\label{eq:metric-energy}
\widetilde{\mathcal E}(f,g) =\sum_{e\in E}\int_{I_e}f'(z)g'(z)\,dz,
\end{equation}
on continuous functions that are absolutely continuous on each 
interval and have square-integrable derivative. For $K\in\mathcal K_{<\infty}$, let $f_K$ minimize the energy $\widetilde{\mathcal E}(f_K, f_K)$ subject to $f_K|_K=1$ and $\sum_{x\in V} {\rm deg}(x) f_K(x)=0$. Following \cite[Definition 3.2 and Lemma 3.4]{goswami2026crifinite}, write
\begin{equation}\label{eq:capacity}
 r_K = \widetilde{\mathcal E}(f_K,f_K) = q \nu_K
\end{equation}
where $q = 2|E| = \sum_{x \in V} {\rm deg}(x)$. This is the {\em 
zero-average} capacity of $K$. It is increasing under set inclusion 
and satisfies
\begin{equation}\label{eq:capacity-lower}
r_K \ge c |K \cap V|, \qquad K\in\mathcal K_{<\infty};
\end{equation}
see \cite[Proposition 3.5]{goswami2026crifinite}.

Let $\widetilde X$ be the canonical diffusion associated with
\eqref{eq:metric-energy}, with law $\mathbf P_y$ when started at 
$y$, and let $H_K$ be its first hitting time of $K$ (see \cite[\S2.3]{goswami2026crifinite}). The Dirichlet martingale functional from 
\cite[Definition 4.1 and Lemma 4.2]{goswami2026crifinite}
has the following explicit expression:
\begin{equation}\label{eq:martingale-functional}
 M_K \stackrel{{\rm def.}}{=} \nu_K\sum_{y\in V} {\rm deg}(y) \mathbf E_y[\varphi(\widetilde X_{H_K})]
     = M_K^{\rm blk} + M_K^{\rm bdr},
\end{equation}
where
\begin{align}
 M_K^{\rm blk}&\stackrel{{\rm def.}}{=}\nu_K\sum_{x\in K\cap V} {\rm deg}(x) \varphi(x),
 \label{eq:bulk-functional}\\
 M_K^{\rm bdr}&\stackrel{{\rm def.}}{=}\nu_K\sum_{y\in V\setminus K} {\rm deg}(y) \mathbf E_y[\varphi(\widetilde X_{H_K})].
 \label{eq:boundary-functional}
\end{align}
In particular, the boundary term is a {\em nonnegative} linear combination of the values of $\varphi$ on $\partial K$. The recurrence of the diffusion on this finite metric graph implies 
$\mathbf P_y[H_K<\infty]=1$. We recall the continuous cluster exploration and its martingale property at the point where they are used in \S\ref{subsec:cls_bm}.

\section{Auxiliary estimates}\label{sec:auxiliary}
We first record an orthogonal decomposition of the 
field that will play a key role in our sprinkling argument. Similar decompositions have been used previously in the context of level-set percolation; see \cite{DCGRS20, MR4642820}. A more 
general version of the following lemma was proved in \cite[Proposition~3.1]{MR4642820}, to which we refer the interested reader for the proof.
%The centering in this decomposition is with respect to the degree measure, as is the zero- average condition for $\varphi$.

\begin{lemma}\label{lem:independent-decomposition}
On an extension of the underlying probability space, there exist independent centered Gaussian 
variables $(Z_x:x\in V)$ with variances $\Cl{C:var}({\rm deg}(x))^{-1}$ for some absolute constant $\Cr{C:var} \in (0, \infty)$, and an independent centered Gaussian field 
$(\eta(x):x\in V)$ such that
\begin{equation}\label{eq:independent-decomposition}
 \varphi(x)=\eta(x)+Z_x-\overline Z,\qquad
 \overline Z=\sum_{x\in V}\pi_x Z_x \:\: \mbox{with} \:\: \pi_x=\frac{{\rm deg}(x)}{q}
\end{equation}
(recall from \eqref{eq:capacity} that $q = 2|E|$). 
\end{lemma}
% \begin{proof}
% Let $L$ be the combinatorial graph Laplacian, $L^\dagger$ its
% Moore--Penrose inverse, and $\Pi=I-\mathbf1\pi^{\mathsf T}$.
% The covariance matrix of $\varphi$ is
% $G=\Pi L^\dagger\Pi^{\mathsf T}$: indeed, the usual free field modulo
% constants has precision $L$, and $\Pi$ chooses its representative with
% degree-weighted average zero. Since the largest eigenvalue of $L$ is at
% most $2d$, and $\Pi\mathbf1=0$, we have, in the order of symmetric matrices,
% \[
%  G\ge\frac1{2d}\Pi\Pi^{\mathsf T}.
% \]
% Consequently $G-\sigma^2\Pi\Pi^{\mathsf T}$ is positive semidefinite.
% Take $\eta$ with this covariance independently of $Z$. Then
% $\eta+\Pi Z$ has covariance $G$, proving \eqref{eq:independent-decomposition}.
% Both covariance matrices annihilate the degree vector, which proves the
% asserted zero-average property of $\eta$.
% Also, $\operatorname{Var}(\overline Z)=\sigma^2\sum_x\pi_x^2\le C/n$, since
% $n\le q\le dn$ and $1\le {\rm deg}(x)\le d$. The Gaussian tail bound gives
% \eqref{eq:mean-noise-tail}. Finally, an additive constant shifts the linear
% interpolation of the endpoint values by the same constant and leaves each
% centered bridge unchanged. This proves \eqref{eq:perturbed-field}.
% \end{proof}
We continue to use $\P$ for the probability measure on the 
extended space supporting these new random fields. Using the 
bound on the variances of $Z_x$ and tail bounds for Gaussian variables, we obtain the following useful estimate.
\begin{corollary}\label{cor:Zbar}
For any $u > 0$, we have
\begin{equation}\label{eq:mean-noise-tail}
 \P[|\overline Z|>u] \le 2\exp(-c |V| u^2).
\end{equation}
\end{corollary}

Notice that, since the field $\varphi$ satisfies the zero-average condition, we also have
$$\sum_{x \in V} {\rm deg}(x) \eta(x)=0$$ almost surely. The field $\varphi_{01}=\eta+Z$ can be extended to the metric graph 
using the same centered Brownian bridges as $\varphi$, so that
\begin{equation}\label{eq:perturbed-field}
 \varphi =  \varphi_{01} - \overline Z \quad\hbox{on }\widetilde{\mathcal G}.
\end{equation}

\smallskip

Next, we bound the `bulk part' of the Dirichlet martingale functional $M_K$, i.e., $M_K^{\rm blk}$ from \eqref{eq:bulk-functional} (cf.~Lemma~7.9 in \cite{goswami2026crifinite}). 
This will be used in \S\ref{subsec:cls_bm} to show 
that $M_{\mathcal C}$, for any {\em small} component $\mathcal C$ of $\widetilde{\mathcal G}^{\ge h}$, is dominated by the 
values of the field $\varphi$ along its boundary $\partial \mathcal C$. %The distinction between the volume cutoff and the error in the following estimate is useful: the latter is fixed before the former is chosen.
\begin{lemma}\label{lem:bulk-control}
For every $b > 0$ and $\varepsilon\in(0,1)$, there exists
$\Cl[c]{c:hep} = \Cr{c:hep}(b, \varepsilon)\in (0, 1)$ such 
that probability of the event
\begin{equation}\label{eq:bulk-event}
 F_0 = F_{0}(b, \varepsilon) \stackrel{{\rm def.}}{=} \Big\{\Big|\sum_{x\in S} {\rm deg}(x) \varphi(x)\Big|
       \le\varepsilon b q\text{ for every }S \subset V
       \text{ with }|S|\le \Cr{c:hep} |V| \Big\}
\end{equation}
is bounded as 
$$\P[F_0^c]\le C e^{-|V|}.$$ In particular, on the 
event $F_0$, we have
\begin{equation}\label{eq:bulk-bound}
 |M_K^{\rm blk}|\le\varepsilon b r_K \text{ whenever $K \in \mathcal K_{< \infty}$ satisfies $|K \cap V| \le \Cr{c:hep} |V|$}.
\end{equation}
%(recall $M_K^{\rm blk}$ from \eqref{eq:bulk-functional}). 
\end{lemma}

\begin{proof}
The variance estimate \eqref{eq:weighted-sum-variance} gives ${\rm Var}\big[\sum_{x\in S}{\rm deg}(x)\varphi(x)\big]\le C|S|$. Thus, for every nonempty $S$ with $|S|\le a |V|$,
\begin{equation*}   
 \P\Big[\Big|\sum_{x\in S}{\rm deg}(x)\varphi(x)\Big|>
                 \varepsilon b q\Big]
 \le2\exp\left(-\frac{c\varepsilon^2 b^2 |V|}{a}\right).
\end{equation*}
For $a\le1/2$, the number of subsets of $V$ of size at most $a|V|$ is at most $\exp(|V|H(a))$, where $H(a)=-a\log a-(1-a)\log(1-a)$. A union bound therefore gives
\begin{equation*}
\P[F_0^c]\le C \exp\left(|V| H(a)-\frac{c\varepsilon^2 b^2|V|}{a}\right).    
\end{equation*}
We can make the exponent smaller than $-|V|$ by choosing $a$ 
small enough depending on $\varepsilon$ and $b$ (in addition to 
$D$ and $\lambda$). The definitions of $M_K^{\rm blk}$ and 
$r_K$ from \eqref{eq:bulk-functional} and \eqref{eq:capacity} 
respectively then yield \eqref{eq:bulk-bound}. 
\end{proof}

Our next result provides a tail bound on the difference in the number of open edges between two nearby levels. 
One useful consequence of this bound is that any {\em large} cluster in $\widetilde{\mathcal G}^{\ge h}$ retains a large sub-cluster in $\widetilde{\mathcal G}^{\ge h + \delta}$ for sufficiently small $\delta$ with high probability (see Lemma~\ref{lem:preservation-large-cluster} in the next section).
\begin{lemma}\label{lem:edge-deletion}
Fix a compact interval $J\subset\mathbb R$, $b>0$ and 
$A > 0$. Put $k=\lceil b|V| /\log |V|\rceil$, $u = \lceil k/(2D)\rceil$ and $s=\lceil k/2\rceil$. For all $h\in J$, $\delta \in (0, 1)$ and $|V|$ large enough depending only on $J, D$ and $\lambda$, we have
\begin{equation}\label{eq:edge-deletion}
 \P[|E^{\ge h}\setminus E^{\ge h + \delta}|\ge k]
 \le\left(\frac{C \delta |V|}{u}\right)^u +
      \left(\frac{C \delta |V| \sqrt{\log |V|}}{s}\right)^s + |V|^{-A}.
\end{equation}
In particular, for $\delta = \kappa(\log |V|)^{-3/2}$ 
and sufficiently small $\kappa = \kappa(b, D, \lambda, J, A)>0$, the left-hand side is at most $2 |V|^{-A}$ 
provided $|V|$ is large enough depending on $J, D, \lambda$ and $b$.
\end{lemma}

The proof of this result rests on the following two 
lemmas.
\begin{lemma}\label{lem:vertex-interval}
For $S\subset V$ with $|S|\le |V|/2$ and real intervals 
$(J_x:x\in S)$ of length at most $\ell$,
\begin{equation}\label{eq:vertex-interval}
 \P[\varphi(x)\in J_x\text{ for every }x\in S]\le(C \ell)^{|S|}.
\end{equation}
Consequently, for $1\le k\le |V|/2$,
\begin{equation}\label{eq:vertex-interval-count}
 \P\left[\big|\{x\in V:\varphi(x)\in J\}\big|\ge k\right]
 \le\left(\frac{C\ell |V|}{k}\right)^k
\end{equation}
for every interval $J$ of length $\ell$.
\end{lemma}

\begin{proof}
In view of the orthogonal decomposition in \eqref{eq:independent-decomposition}, we may first condition on 
the field $\eta$; then \eqref{eq:vertex-interval} follows from 
the bound below, holding uniformly over all intervals $(J_x : x \in S)$ of maximum length $\ell$:
\begin{equation}\label{eq:ZxbarZ}
\P[Z_x - \overline Z \in J_x\text{ for every }x\in S]\le(C \ell)^{|S|}.
\end{equation}
To this end notice that, since $|S| \le |V| / 2$ and $(Z_x : x \in V)$ are distributed as independent $N(0, \Cr{C:var}({\rm deg}(x))^{-1})$ variables by Lemma~\ref{lem:independent-decomposition}, the conditional distribution of $Z_x$ given 
$\mathbf Z_{x, S} = ((Z_y : y \in S \setminus \{x\}), \overline Z)$ is $N(\mu, \sigma_{x, S}^2)$ for some $\mu = \mu(\mathbf Z_{x, S})$ and $\sigma_{x, S} \ge c \in (0, \infty)$. In 
particular, the conditional density of $Z_x$ is bounded above by some $C \in (0, \infty)$ and hence
\begin{equation*}
\P[Z_x - \overline Z \in J_x \, | \, \mathbf Z_{x, S} ] \le C |J_x| \le C \ell
\end{equation*}
for all $x \in S$ which yields \eqref{eq:ZxbarZ}. 
\eqref{eq:vertex-interval-count} follows from \eqref{eq:vertex-interval} via a union bound over the $\binom{|V|}{k}\le(e|V|/k)^k$ possible sets $S$.
\end{proof}

\begin{lemma}\label{lem:edge-interval}
Fix a compact interval $J\subset\mathbb R$ and $A>0$. For $h\in J$ and $\delta \in (0, 1)$, let $N(h, \delta)$ be the number of edges $e=\{x,y\}$ for which $\varphi(x),\varphi(y)\ge h + \delta$ and $\min_{z\in\overline I_e}\varphi(z)\in[h,h + \delta)$. Then, for $1\le k\le|E|$ and $|V|$ large enough 
depending only on $J, D$ and $\lambda$, we have
\begin{equation}\label{eq:edge-interval}
 \P[N(h, \delta)\ge k]
 \le\left(\frac{C(A)\, \delta |V|\sqrt{\log |V|}}{k}\right)^k + |V|^{-A}.
\end{equation}
\end{lemma}

\begin{proof}
Recall that conditionally on the field on $V$, the values on 
different edges are distributed as independent standard 
Brownian bridges joining their respective endpoint values. For endpoint values $a,b\ge h+\delta$, the reflection principle for the Brownian bridge (see, e.g., \cite[Chap.IV, \S26, p.~67]{BorodinBrown}) gives the conditional probability
\begin{equation}\label{eq:cond_prob}
e^{-2(a-h-\delta)(b-h-\delta)}-e^{-2(a-h)(b-h)} \le 2\delta(a+b-2h).    
\end{equation}
The inequality follows simply by factoring out the first exponential and using $1-e^{-t}\le t$ with $t=2\delta(a+b-2h)-2\delta^2$. Now let $F_R \stackrel{{\rm def.}}{=} \{\max_{x\in V}|\varphi(x)|\le R\sqrt{\log |V|}\}$ for any $R \ge 0$. The uniform variance bound for each $\varphi(x)$ (cf.~\eqref{eq:weighted-sum-variance}) together with standard 
Gaussian tail estimates yield the upper bound
\begin{equation}\label{eq:FR}
\P[F_R^c] \le |V|^{-A}    
\end{equation}
for large enough $R=R(A)$. On the event $F_R$, the conditional probability in \eqref{eq:cond_prob} is at most $$\big(R \sqrt{\log |V|} + C(J)\big) \, \delta.$$ Conditional independence and a union bound over $\binom{|E|}{k}\le(e D |V|/ k)^k$ sets of edges together with \eqref{eq:FR} therefore give 
\eqref{eq:edge-interval}.
\end{proof}
We can now give the
\begin{proof}[Proof of Lemma~\ref{lem:edge-deletion}]
%Any edge $e = \{x, y\}$ in $E^{\ge h} \setminus E^{\ge h + \delta}$ either has an endpoint with the value of $\varphi$ lying in $[h, h + \delta)$ or is counted by $N(h, \delta)$. Thus, almost surely,
Notice that, by the definitions of $N(h, \delta)$ and 
$E^{\ge h}$, we have the bound:
\begin{equation*}
|E^{\ge h}\setminus E^{\ge h + \delta}|
 \le D\big|\{x:\varphi(x)\in[h, h + \delta)\}\big| + N(h, \delta).   
\end{equation*}
If the left-hand side is at least $k$, either the first 
term on the right-hand side is bounded below by $u$ or 
$N(h, \delta)\ge s$. Adding up the bounds given by 
Lemmas~\ref{lem:vertex-interval} and~\ref{lem:edge-interval}, we obtain \eqref{eq:edge-deletion}.
\end{proof}

We will use the following elementary observation to 
transfer large clusters between nearby levels.

\begin{lemma}\label{lem:deletion-components}
Consider two graphs $H'=(V,F')$ and $H=(V,F)$ with $F'\subset F$. Suppose that every connected component 
of $H'$ has either at most $L$ vertices or at least $a'|V|$ vertices. If $|F\setminus F'|+1<a|V|/L$, then every 
connected component of $H$ with at least $a|V|$ vertices 
contains a component of $H'$ with at least $a'|V|$ 
vertices.
\end{lemma}

\begin{proof}
Deleting $q'$ edges from a connected graph leaves at most $q'+1$ connected components. If a connected component of $H$ with at least $a|V|$ vertices contained
no connected component of $H'$ with at least $a'|V|$ vertices, all of such components would have at most $L$ vertices. Its volume would consequently be at most 
$(|F\setminus F'|+1)L<a|V|$, which is a contradiction.
\end{proof}

\section{The supercritical phase}\label{sec:supcri}
In this section we prove Theorem~\ref{thm:uniq_giant} by combining two main ingredients: Proposition~\ref{prop:no_intermediate_uniform} and 
Theorem~\ref{thm:ubd_F2}, established in \S\ref{subsec:cls_bm} 
and \S\ref{subsec:uniq_giant} respectively. The former rules out the existence of any cluster of `intermediate' volume, 
i.e., volume in the range $[A \log |V|, a |V|]$ for positive 
constants $a = a(h, D, \lambda)$ and $A = A(h, D, \lambda)$. 
The complementary result, Theorem~\ref{thm:ubd_F2}, shows that 
with high probability there can be at most one cluster (the 
giant cluster) with volume at least $a|V|$, using this {\em 
same} choice of constant $a$.

\subsection{Non-existence of intermediate 
clusters}\label{subsec:cls_bm}
In this subsection, we prove the absence of clusters between logarithmic and 
small linear volume, uniformly over a compact interval of 
supercritical levels. This will also allow us to apply the 
result at the slightly higher level used in the sprinkling 
argument. %Throughout this subsection, we write $q = 2|E|$ (recall \eqref{eq:capacity}) and $b=|h|/2$, where $h<0$ is fixed.
\begin{lemma}\label{lem:brownian_last_crossing}
Let $(B_t)_{t\ge0}$ be a standard Brownian motion under 
$\mathbf P^B_0$. For $\beta>0$ and $s>0$,
\begin{equation}\label{eq:brownian_last_crossing}
 \mathbf P^B_0[\text{there exists }t\ge s\text{ such that }B_t\le-\beta t]
 =2\Phi(-\beta\sqrt{s})\le \exp(-\beta^2s/2),
\end{equation}
where $\Phi$ denotes the standard normal distribution function.
\end{lemma}
\begin{proof}
This follows from standard properties of Brownian motion; see, e.g., the proof of Lemma~7.3 in \cite{goswami2026crifinite}. We 
omit the details for brevity.
% Set $Y_t=B_t+\beta t$. For $y>0$, Brownian motion with drift $\beta$, started at $y$, hits zero with probability $e^{-2\beta y}$. To recall a proof, $e^{-2\beta Y_t}$ is a martingale, and optional stopping at the exit time from
% $(0,R)$ gives the probability of hitting zero before $R$ as
% $(e^{-2\beta y}-e^{-2\beta R})/(1-e^{-2\beta R})$; let $R\to\infty$.
% The Markov property at time $s$ now yields
% \[
%  \mathbf P^B_0[\inf_{t\ge s}Y_t\le0]
%  =\Phi(-\beta\sqrt{s})+
%  \int_0^\infty e^{-2\beta y}
%        \frac{\exp(-(y-\beta s)^2/(2s))}{\sqrt{2\pi s}}\,dy
%  =2\Phi(-\beta\sqrt{s}).
% \]
% The usual Gaussian tail bound gives the inequality.
\end{proof}

Let $h < 0$ be fixed and $b=|h|/2$. Consider the event $$F_0' = F_{0}(b, \tfrac18)$$
introduced in \eqref{eq:bulk-event}. Applying 
Lemma~\ref{lem:bulk-control},
% with $\varepsilon = 1/16$ and decrease its volume cutoff if necessary to choose
% $a\in(0,\min\{1/2,1/(4d)\})$. Denote the resulting event by
% \begin{equation}\label{eq:bulk_good_uniform}
%  \mathcal B_h = \left\{\left|\sum_{x\in S}{\rm deg}(x)\varphi(x)\right|
%  \le bq/8\text{ for every }S\subset V\text{ with }|S|\le an\right\}.
% \end{equation}
we obtain 
\begin{equation}\label{eq:bnd_F0'}
\P[(F_0')^c]\le C e^{-|V|}.    
\end{equation}
Further, on the event $F_0'$, we have $|M_K^{\rm blk}|\le br_K/8$ whenever $|K \cap V| \le a|V|$ for $$a = \min\{1/2,1/(4D), \Cr{c:hep}(b, \tfrac 1{16})\}$$ 
($\Cr{c:hep}$ is as in Lemma~\ref{lem:bulk-control} and $M_K^{\rm blk}$ is from \eqref{eq:bulk-functional}). In the 
remainder of this Section, we work with this particular choice 
of $a$.

\begin{proposition}\label{prop:no_intermediate_uniform}
There exists $A = A(h,D,\lambda) \in (0, \infty)$ such that, 
for every $h' \in [h,h/2]$,
\begin{equation}\label{eq:no_intermediate_uniform}
 \P\big[\text{there exists a cluster of }\mathcal G^{\ge h'}
       \text{ with volume in }[A\log |V|,a |V|]\big]\le C(h) \,  |V|^{-2}.
\end{equation}
\end{proposition}
\begin{proof}
We first recall that the extension of $\varphi$ to an edge is 
its linear interpolation between the endpoint values plus an 
independent standard Brownian bridge. The variances of 
$\varphi(x)$, for $x \in V$, are bounded in view of 
\eqref{eq:weighted-sum-variance}. A union bound over 
vertices and edges, together with Gaussian and
Brownian bridge tail bounds, therefore gives
\begin{equation}\label{eq:metric_max_tail}
 \P \big[ \mbox{$\sup_{z\in\widetilde{\mathcal G}}$} \:\, \varphi(z)> L \big]
       \le C |V| e^{-cL^2},\qquad L\ge1.
\end{equation}
For sufficiently large $|V|$ depending only on $h$ and $D$, so that $q = 2|E| \ge32 D$ and $bq/(16D)\ge1$, let
$$\mathcal H = \mathcal H(h) \stackrel{{\rm def.}}{=} \big\{\mbox{$\sup_{z\in\widetilde{\mathcal G}}$} \:\, \varphi(z)\le bq/(16 D)\big\}.$$
\eqref{eq:metric_max_tail} then implies
$$\P[\mathcal H^c]\le C |V| e^{-c(h) |V|^2}.$$

\medskip

Now fix $x\in V$, an integer $N>1$, and $h' \in[h,h/2]$. We will use the exploration $(K_t)_{t\ge0}=(K_t(x,h',N))_{t\ge0}$ 
described in \cite[Proposition 6.1]{goswami2026crifinite} with
$K_0=\partial B(x,1/N)$, and put $\tau=\tau_{1,1}$.
The sets $K_t$ are non-decreasing and continuous in the Hausdorff topology,
are contained in $\widetilde{\mathcal G}\setminus B(x,1/N)$, and each
component of $K_t$ intersects $K_0$.  Their measurability property makes
$\widetilde M_t=M_{K_t}$ (see \eqref{eq:martingale-functional}) a continuous martingale in the exploration
filtration.  More precisely, by
\cite[Lemma 6.3, Theorem 6.4]{goswami2026crifinite},
\begin{equation}\label{eq:exploration_clock_recall}
 \langle\widetilde M\rangle_t=r_{K_t}-r_{K_0},
 \qquad
 (\widetilde M_{T_s})_{s\ge0}
 \stackrel{\mathrm{law}}=
 (B_{r_{K_0}+s\wedge\langle\widetilde M\rangle_\infty})_{s\ge0},
 \quad
 T_s=\inf\{t:\langle\widetilde M\rangle_t>s\},
\end{equation}
where $r_K$ is as in \eqref{eq:capacity} and $(B_t)_{t \ge 0}$ is a standard Brownian motion. %The equality is understood with the continuous terminal value when the quadratic variation stops.
Let $F_x(N) = F_x(h', N)$ be the event that the field stays on 
the same side of $h'$ as $\varphi(x)$ throughout 
$\overline{B(x,1/N)}$. %: it is at least $u$ there if $\varphi(x)\ge u$, and strictly less than $u$ there otherwise. 
This is the event used in
\cite[(6.9), Lemma 6.2]{goswami2026crifinite}.
Continuity of the field and the nondegenerate Gaussian law of
$\varphi(x)$ imply $\P[F_x(N)^c]\to0$ as $N\to\infty$.
On this event, if
$2\le k\le|\mathscr C_x^{\ge h'}|\le a|V|$, then
\begin{equation}\label{eq:terminal_volume}
 k-1\le |K_\tau \cap V| = |\mathscr C_x^{\ge h'}|-1\le a|V|.
\end{equation}
Indeed, by \cite[Proposition 6.1(3-i)]{goswami2026crifinite}, $K_\tau$ is the union of $K_0$ and the complete 
components at level $h'$ which meet $K_0$ outside $B(x,1/N)$; 
on the event $F_x(N)$ these form precisely the component of $x$ 
with the vertex $x$ removed.

\smallskip

Now set $K=K_\tau$ and $D_K=\sum_{y\in K\cap V}{\rm deg}(y)$. 
Is follows from the special boundary property in
\cite[Proposition 6.1(3-ii)]{goswami2026crifinite} that
$$\partial K\setminus\{\varphi\le h' \}\subset K_0,
 \qquad
 \mathbf P_y[\widetilde X_{H_K}\in K_0\setminus\{\varphi\le h'\}]=0
 \quad(y\in V\setminus\{x\}).$$
Thus the contribution from points $z \in \partial K$ with $\varphi(z) \ge h'$ to the expression for $M_K^{{\rm bdr}}$ in \cite[(4.3), Lemma 4.2]{goswami2026crifinite} (see also \eqref{eq:boundary-functional}) can come only from a walk started at $x$. This gives the bound
\begin{equation}\label{eq:terminal_boundary_upper}
 M_K^{\mathrm{bdr}}
 \le\nu_K\Big\{h'(q-D_K)
              + {\rm deg}(x)\big(\sup_{z\in K_0}\varphi(z)-h'\big)_+\Big\}.
\end{equation}
On $\mathcal H$, since $|h'|\le 2b$ and $q\ge32D$,
$${\rm deg}(x)\big(\sup_{z\in K_0}\varphi(z)-u\big)_+
 \le bq/16 + 2bD\le bq/8.$$
Moreover, $D_K\le D a |V|\le q/4$ by our choice of $a$. Using 
$h' \le -b$, \eqref{eq:terminal_boundary_upper}, and the bound 
on $|M_K^{{\rm blk}}|$ on the event $F_0'$ as recorded below 
\eqref{eq:bnd_F0'} (recall from \eqref{eq:capacity} that $r_K = q \nu_K$), we obtain
\begin{equation}\label{eq:terminal_martingale_negative}
 M_K=M_K^{\mathrm{blk}}+M_K^{\mathrm{bdr}}
 \le\nu_K\{bq/8-3bq/4+bq/8\}
 =-\frac b2r_K.
\end{equation}

%The quadratic variation and time-change statements in\cite[Lemma 6.3, Theorem 6.4]{goswami2026crifinite}, based on theDubins--Schwarz theorem \cite[Chapter V, Theorem 1.7]{revuz2013continuous}, identify $M_{K_t}$ with a standard Brownian motion at time $r_t=r_{K_t}$(up to its terminal quadratic variation, with the initial Brownian time$r_0=r_{K_0}$ included).  Consequently,
Translated in terms of the Brownian motion $(B_t)_{t \ge 0}$ in 
\eqref{eq:exploration_clock_recall}, 
\eqref{eq:terminal_martingale_negative} implies that, 
this Brownian motion satisfies $B_r\le-br/2$ at $r=r_K$.
By \cite[Proposition 3.5]{goswami2026crifinite} and
\eqref{eq:terminal_volume}, $r\ge \Cl[c]{c:cap}(k-1)$ for a 
constant $\Cr{c:cap} \in (0, \infty)$. Lemma~\ref{lem:brownian_last_crossing} therefore yields
$$\P[k\le|\mathscr C_x^{\ge u}|\le a|V|]
 \le 2e^{-|V|} + C |V| e^{-c(h)|V|^2}
       +e^{-b^2c_0(k-1)/8}+\P[F_x(N)^c].$$
The last term tends to zero as $N\to\infty$ by
\cite[Lemma 6.2]{goswami2026crifinite}. All the other bounds are independent
of $N$ and of $h' \in [h,h/2]$. Sending $N\to\infty$, summing over $x\in V$, and choosing $k=\lceil A\log |V|\rceil$ with $A$ sufficiently large proves
\eqref{eq:no_intermediate_uniform}.
\end{proof}

\subsection{Uniqueness of the giant component}\label{subsec:uniq_giant}
We now prove uniqueness at the same linear scale as in
Proposition~\ref{prop:no_intermediate_uniform}. Throughout this subsection, fix $h<0$ and choose $a=a(h, D, \lambda)$ and $A=A(h, D, \lambda)$ so that the conclusion of that
proposition holds, with these same constants, at every level $h'\in[h,h/2]$.

\begin{theorem}[Uniqueness of the giant component]\label{thm:ubd_F2}
There exists $C(h) \in (0, \infty)$ such that
\begin{equation}\label{eq:ubd_F2}
 \P\big[\text{there are two distinct clusters in $\mathcal G^{\ge h}$
 of volume at least $a|V|$}\big]\le C(h) |V|^{-2}.
\end{equation}
\end{theorem}

Recall the decomposition from Lemma~\ref{lem:independent-decomposition}:
\begin{equation}\label{def:psi_01}
 \varphi=\eta+Z-\overline Z,\qquad
 \varphi_{01}=\eta+Z=\varphi+\overline Z,
 \qquad
 \overline Z=\sum_{x\in V} \pi_x Z_x.
\end{equation}
Here the $Z_x$ are independent $N(0,\sigma_x^2)$ variables with $\sigma_x^2 = \Cr{C:var} \, {{\rm deg}(x)}^{-1}$, independent of $\eta$. As noted above \eqref{eq:perturbed-field}, on each edge we use the same centered Brownian bridge for $\varphi$ and 
$\varphi_{01}$. Consequently, the last equality in \eqref{def:psi_01} holds everywhere
on $\widetilde{\mathcal G}$. Denote by $E_{01}^{\ge h'}$ the open edges for $\varphi_{01}$ at level $h'$, i.e., the analogue 
of $E^{\ge h'}$ with $\varphi_{01}$ in place of $\varphi$ (cf.~\eqref{eq:cable-level-set}).

Let $0<\delta<\min\{1,-h/4\}$. Let $F_1$ be the event that any two
clusters of $(V,E_{01}^{\ge h+\delta})$ with volume at least $a|V|$ are
connected in $(V,E_{01}^{\ge h+\delta/2})$. Let $F_2$ be the event that
every cluster of $\mathcal G^{\ge h}$ with volume at least $a|V|$ contains
a cluster of $(V,E_{01}^{\ge h+\delta})$ with volume at least $a |V|$.
Finally, put $F_3=\{|\overline Z|\le\delta/2\}$. In view of \eqref{def:psi_01}, we have the following inclusions on the 
event $F_3$:
\begin{equation}\label{eq:shift-inclusions}
 \{\varphi\ge h+2\delta\}\subset
 \{\varphi_{01}\ge h+\delta\}\subset\{\varphi\ge h\},
 \qquad
 \{\varphi_{01}\ge h+\delta/2\}\subset\{\varphi\ge h\}.
\end{equation}
It follows that on $F_1\cap F_2\cap F_3$ there is at most one
cluster of $\mathcal G^{\ge h}$ of volume at least $a|V|$. 
Corollary~\ref{cor:Zbar} already provides a bound on $\P[ F_3^c ]$. We now estimate the probabilities of the complements of $F_2$ and $F_1$ separately.

\begin{lemma}\label{lem:preservation-large-cluster}
There exist $c = c(h), C = C(h) \in (0, \infty)$ such that
\begin{equation}\label{eq:preservation-large-cluster}
 \P[F_2^c]\le
 C(C\delta\log^{3/2}|V|)^{c |V|/\log |V|}
       +C|V|^{-2}+C\exp(-c\delta^2 |V|).
\end{equation}
\end{lemma}
\begin{proof}
Let $F'_4$ be the event that every cluster of $\mathcal G^{\ge h}$ of
volume at least $a|V|$ contains a cluster of $\mathcal G^{\ge h+2\delta}$
of volume at least $a|V|$. On $F_3$, the inclusions
\eqref{eq:shift-inclusions} give $F'_4\subset F_2$.
By Proposition~\ref{prop:no_intermediate_uniform}, except on an event of probability at most
$C|V|^{-2}$, every cluster of $\mathcal G^{\ge h+2\delta}$ has volume either
less than $A\log |V|$ or at least $a|V|$. On this event,
Lemma~\ref{lem:deletion-components}, with $L=A\log |V|$ and $a'=a$,
shows that failure of $F'_4$ requires at least $a|V|/(A\log |V|)-1$ edges
in $E^{\ge h}\setminus E^{\ge h+2\delta}$. Thus, for $|V|$ large enough depending only on $h, D$ and $\lambda$, this difference contains at least $b|V|/\log |V|$ edges, 
with $b = a/(2A)$. Lemma~\ref{lem:edge-deletion} bounds the probability of this event by the first two terms in
\eqref{eq:preservation-large-cluster}. Finally, by Corollary~\ref{cor:Zbar},
\begin{equation}\label{eq:shift-tail}
 \P[F_3^c]\le C\exp(-c\delta^2 |V|).
\end{equation}
Putting all these bounds together, we obtain 
\eqref{eq:preservation-large-cluster}.
\end{proof}

We need the following combinatorial lemma for bounding $\P[ F_1^c ]$ which uses the expansion property of the graph $\mathcal G$.
\begin{lemma}\label{lem:many-short-paths}
There exist $\rho=\rho(a, D, \lambda)>0$ and an integer
$\ell=\ell(a, D,\lambda)\ge1$ such that any disjoint vertex sets
$S,T\subset V$ with $|S|,|T|\ge a |V|$ can be joined by at least
$2\rho |V|$ pairwise vertex-disjoint paths of length at most 
$\ell$, for all sufficiently large $|V|$ depending only on $a, 
D$ and $\lambda$.
\end{lemma}
\begin{proof}
By property \eqref{eq:edge-expansion}, $|\partial_e U|\ge c |U|$ whenever $|U|\le |V|/2$. Let $W$ be a vertex set separating $S$ from $T$, where vertices of $S\cup T$ may belong to $W$. If $|W|< a|V|/2$, let $U$ be the union of the components of $V\setminus W$ meeting $S\setminus W$. Then $U$ contains $S\setminus W$ and its complement outside $W$
contains $T\setminus W$. One of these two sets has volume at 
most $|W|/2$, and both have volume at least $a|V|/2$. Its edge 
boundary is contained in the edges incident to $W$, whence
$D|W|\ge c a|V|/2$. Thus every such separator has at least $b_0|V|$
vertices, where $b_0>0$ depends only on $a,D,\lambda$. Menger's theorem (see, e.g., \cite[Chapter 5]{bang2008digraphs}) yields at least $\lfloor b_0 |V|\rfloor$ pairwise vertex-disjoint $S$--$T$ paths. The sum of their numbers of vertices is at most $|V|$. Choosing $\ell$ sufficiently large leaves at least
$b_0 |V|/4$ paths of length at most $\ell$. Reducing $\rho$ if
necessary proves the assertion.
\end{proof}

For each pair $S,T$ as above, choose one such collection of short
paths according to a fixed deterministic ordering. A path is
{\em $L$-good} if $\eta(x)>-L$ at every vertex on the path. We choose
$L=L(h, D,\lambda)<\infty$ so that
\begin{equation}\label{eq:good-background}
 \P\big[|\{x\in V:\eta(x)\le-L\}|>\rho |V|\big]
       \le C\exp(-c|V|).
\end{equation}
Here and below $\rho$ is as in Lemma~\ref{lem:many-short-paths}.
To verify \eqref{eq:good-background}, observe that
${\rm Cov}(\eta)\le {\rm Cov}(\varphi)$ as
quadratic forms. Therefore the variance estimate \eqref{eq:weighted-sum-variance} gives
$${\rm Var}\Big[\sum_{x\in U}{\rm deg}(x)\eta(x)\Big]
       \le C|U|,\qquad U\subset V.$$
If all vertices of a fixed set $U$ satisfy $\eta(x)\le-L$,
then $\sum_{x\in U}{\rm deg}(x)\eta(x)\le-L|U|$. A Gaussian tail bound
and a union bound over subsets of $V$ show that the left-hand side
of \eqref{eq:good-background} is at most
$2^{|V|}\exp(-cL^2\rho |V|)$, up to an adjustment for integer parts.
Taking $L$ large proves the claim. On the complementary event in
\eqref{eq:good-background}, each chosen collection contains at least $\rho |V|$-many $L$-good paths.

\begin{lemma}\label{lem:sprinkling-connection}
There exist $c = c(h), C = C(h) \in (0, \infty)$ and an integer $\ell\ge1$, depending only on $h, D,\lambda$, such that
\begin{equation}\label{eq:sprinkling-connection}
 \P[F_1^c]\le C\exp(-c |V| \delta^{3\ell}) + C\exp(-c|V|).
\end{equation}
\end{lemma}
\begin{proof}
Put $t=h+\delta$, $s=h+\delta/2$, and $r=h+3\delta/4$.
We describe precisely the information on which we condition.
Reveal $\eta$ on $V$, the indicators
$H_x=1_{\{\varphi_{01}(x)\ge t\}}$, the field values
$\varphi_{01}(x)$ at vertices with $H_x=1$, and all indicators
$Q_e=1_{\{e\in E_{01}^{\ge t}\}}$. Denote the resulting
$\sigma$-algebra by $\mathcal H$.

Conditionally on these random variables, the values at vertices 
with $H_x=0$ are independent, and at such a vertex $x$ the law is that
of $\eta(x)+Z_x$ conditioned to be less than $t$.
Indeed, an edge incident to such a vertex is necessarily closed
at level $t$, so its indicator $Q_e$ supplies no additional
information. At vertices with $H_x=1$ the field values have been
fixed. The centered bridges are independent conditionally on their
endpoints and the indicators $Q_e$; only bridges with two endpoints
above $t$ are subject to a nontrivial restriction. These assertions
follow directly from the independent Gaussian variables and
independent bridge construction.

At a vertex $x$ with $H_x=0$ and $\eta(x)>-L$, we have the 
uniform lower bound
\begin{equation}\label{eq:truncated-vertex-sprinkling}
 \P[\varphi_{01}(x)\ge r\mid\mathcal H]
       \ge c\delta\exp\{-C(L+|h|+1)^2\}.
\end{equation}
For completeness, we first condition on the field $\eta$ which 
converts the problem to estimating conditional probabilities 
of the form $\P[Z_x \ge q - w \mid Z_x < q]$. Writing $\Phi$ 
for the standard normal distribution function, note that this probability is decreasing in $q$ for fixed $w>0$. This follows by differentiating its expression
$1-\Phi((q-w)/\sigma_x)/\Phi(q/\sigma_x)$, using that
$\Phi'/\Phi$ is decreasing. Apply this with $w=\delta/4$ and
$q=t-\eta(x)\le L+|h|+1$, and integrate the Gaussian density
over the interval of length $w$ at that upper endpoint. This
proves \eqref{eq:truncated-vertex-sprinkling}.

We also need a lower bound for the bridges after conditioning on
their level-$t$ status. If an edge $e = \{x, y\}$ satisfies 
$\varphi(x) = u, \varphi(y) = v$ with $u,v\ge t$ and $Q_e=0$, 
the expressions for the tail probabilities of the minimum of a 
Brownian bridge gives (recall \eqref{eq:cond_prob}):
\begin{align}\label{eq:conditional-bridge-sprinkling}
 &\P[\min_{I_e}\varphi_{01}\ge s\mid
              \varphi_{01}(x)=u,\varphi_{01}(y)=v,Q_e=0]\notag\\
 &\hspace{1cm}=1-\exp\big(-2[(u-s)(v-s)-(u-t)(v-t)] \,\big)
       \ge1-\exp\big(-2(t-s)^2\big).
\end{align}
If $Q_e=1$, the edge is already open at level $s$.
If at least one endpoint is below $t$, the condition $Q_e=0$ is
automatic and the bridge retains its unconditioned law. Once its
endpoints are at least $r$, its probability of staying above $s$
is at least $1-\exp(-2(r-s)^2)$. Thus in all cases, conditionally on endpoint values at least $r$, the required bridge event has
probability at least $c\delta^2$.

The clusters at level $t$ and the chosen collections of paths
between them are $\mathcal H$-measurable, as is the collection of
$L$-good paths. The endpoints of such a path belong to clusters
of volume at least $a|V|>1$, so their field values are already at
least $t$. For an $L$-good path of at most $\ell$ edges, we use
\eqref{eq:truncated-vertex-sprinkling} for the internal vertices
below $t$, and then the bridge bound just proved for its edges.
Therefore, the conditional probability that all its vertices 
are at least $r$ and that it is open at level $s$ is at least
\begin{equation}\label{eq:one-path-sprinkling}
 p_\delta=(c\delta^3)^{\ell},
\end{equation}
where the constant $c$ may depend on $h, D,\lambda$.

For vertex-disjoint paths these events are independent conditionally
on $\mathcal H$. To see this explicitly, all high vertex values
have been fixed, the remaining vertex values are independent
truncated Gaussians, and each bridge used by a path belongs only
to that path. Conditions on edges outside the chosen paths
introduce no dependence between their low vertex values, since
any such edge incident to a low vertex has $Q_e=0$ automatically.
There are at most $1/a$ clusters of volume at least $a|V|$.
On the event in \eqref{eq:good-background} that there are at most $\rho |V|$ `bad' vertices, the probability that a specified pair of these clusters is not joined by one of its good paths is at most
$$(1-p_\delta)^{\rho |V|}\le\exp(-\rho |V| p_\delta).$$
A conditional union bound over at most $a^{-2}$ pairs, 
followed by \eqref{eq:good-background}, then proves 
\eqref{eq:sprinkling-connection}.
\end{proof}

\begin{proof}[Proof of Theorem~\ref{thm:ubd_F2}]
For all sufficiently large $|V|$ depending only on $h, D$ and $\lambda$, choose
$$\delta=\kappa(\log |V|)^{-3/2},$$
where $\kappa>0$ is small enough that the constant $C$ in
\eqref{eq:preservation-large-cluster} satisfies $C\kappa\le1/2$. Increasing the lower bound on $|V|$ ensures
$\delta<\min\{1,-h/4\}$. By \eqref{eq:shift-inclusions}, the event in \eqref{eq:ubd_F2} is contained in
$F_1^c\cup F_2^c\cup F_3^c$. Lemmas
\ref{lem:preservation-large-cluster} and
\ref{lem:sprinkling-connection}, together with
\eqref{eq:shift-tail}, bound its probability by
$$C|V|^{-2}+C\exp(-c |V|/\log |V|)
 +C\exp(-c|V|/\log^3|V|)
 +C\exp\{-c|V|/(\log |V|)^{9\ell/2}\}$$
 with implicit dependence on $h$ in each of the constants in addition to $D$ and $\lambda$ which is precisely the theorem.
\end{proof}

Finally we give the
\begin{proof}[Proof of Theorem~\ref{thm:uniq_giant}]
Applying \cite[Theorem 1.2, (1.4)]{goswami2026crifinite} at level $h < 0$, we get
\begin{equation*}
\P[\, |\mathscr C_{\max}^{\ge h}| < \Cl[c]{c:giant} \min(|h|, 1)  \, |V| \,]
 \le C \min(|h|, 1)^{-1/5}|V|^{-1/15}  
\end{equation*}
for some $\Cr{c:giant}, C \in (0, \infty)$. Without any loss of generality, we can choose $a$ in 
Proposition~\ref{prop:no_intermediate_uniform} so that $a\le \Cr{c:giant} \min(|h|, 1)$. Theorem~\ref{thm:uniq_giant} now follows from the previous display combined with Proposition~\ref{prop:no_intermediate_uniform} and Theorem~\ref{thm:ubd_F2}.
\end{proof}

\noindent \textbf{Acknowledgments.} SG's research is supported by the ANRF grant \sloppy
ANRF/ARGM/2025/001947/MTR and a grant from the Department of Atomic Energy, Government of India, under Project Identification No.~RTI4014. We thank Debapratim Banerjee, Soumendu Sundar Mukherjee, Guillaume Conchon-Kerjan, Pierre-Fran\c{c}ois Rodriguez and Gordon Slade for insightful 
discussions at various stages of this project. We are also grateful to Rajat Subhra Hazra and Remco van der Hofstadfor their valuable comments. SG acknowledges the hospitality of the Statistics and Mathematics Unit, Indian Statistical Institute (ISI), Kolkata on 
several occasions. DP thanks the School of Mathematics, Tata Institute of 
Fundamental Research (TIFR), Mumbai for their generous hospitality.

\bigskip

\noindent \textbf{AI use Statement.} No mathematical content was generated 
or found with the assistance of artificial intelligence.

\end{document}